\documentclass[a4paper,english,oneside]{article}
\usepackage[T1]{fontenc}
\usepackage[utf8]{inputenc} 

\usepackage{geometry}
\usepackage{amsmath}
\usepackage{amsfonts}
\usepackage{amssymb}
\usepackage{amsthm}
\usepackage[english]{babel}
\usepackage{color}
\usepackage{lmodern}
\usepackage{pstricks-add}
\usepackage{hyperref}
\hypersetup{hidelinks}

\usepackage{mathscinet}
\usepackage{enumitem}
\usepackage[cal=boondoxo,bb=ams]{mathalfa}
\usepackage{microtype}

\usepackage{mathtools}
\usepackage{esint}

\theoremstyle{plain}
\newtheorem{theorem}{Theorem}[section]
\newtheorem{lemma}[theorem]{Lemma}

\newtheorem{proposition}[theorem]{Proposition}

\theoremstyle{definition}

\theoremstyle{remark}
\newtheorem{remark}{Remark}

    \DeclareMathOperator\supp{supp}
    \DeclareMathOperator\meas{meas}

    \DeclareMathOperator\Str{Str}
    \DeclareMathOperator\Gla{Gla}
    
    \DeclareMathOperator\Sob{Sob}

\begin{document}

\title{A blow-up result for a generalized Tricomi equation with nonlinearity of derivative type}

\author{Sandra Lucente$^\mathrm{a}$, Alessandro Palmieri$^\mathrm{b}$
}

\date{\small{$^\mathrm{a}$ Department of Physics, University of Bari, Via E. Orabona 4, 70125 Bari, Italy}\\ \small{$^\mathrm{b}$ Department of Mathematics, University of Pisa,  Largo B. Pontecorvo 5, 56127 Pisa, Italy}
}

\maketitle

\begin{abstract}

In this note, we prove a blow-up result for a semilinear generalized Tricomi equation with nonlinear term of derivative type, i.e., for the equation $\mathcal{T}_\ell u=|\partial_t u|^p$, where $\mathcal{T_\ell}=\partial_t^2-t^{2\ell}\Delta$.
 Smooth solutions blow up in finite time for positive Cauchy data when the exponent $p$ of the nonlinear term is below $\frac{\mathcal{Q}}{\mathcal{Q}-2}$, where $\mathcal{Q}=(\ell+1)n+1$ is the quasi-homogeneous dimension of the generalized Tricomi operator $\mathcal{T}_\ell$.
Furthermore, we get also an upper bound estimate for the lifespan.


\end{abstract}

\begin{flushleft}
\end{flushleft}

\begin{flushleft}
\textbf{AMS Classification (2010)} Primary:  35B44, 35L71; Secondary: 35B33, 35C15
\end{flushleft}

\section{Introduction}

In the present work, we prove a blow-up result for the semilinear generalized Tricomi equation with nonlinearity of derivative type, namely,
\begin{align}\label{semilinear CP derivative type}
\begin{cases} \partial_t^2 u-t^{2\ell}\Delta u =|\partial_t u|^p, &  x\in \mathbb{R}^n, \ t>0,\\
u(0,x)=\varepsilon u_0(x), & x\in \mathbb{R}^n, \\ \partial_t u(0,x)=\varepsilon u_1(x), & x\in \mathbb{R}^n,
\end{cases}
\end{align} where $\ell> 0$, $p>1$ and $\varepsilon$ is a positive constant describing the size of Cauchy data.

The semilinear wave equation 
\begin{align}\label{semilinear wave f}
\begin{cases} \partial_t^2 u-\Delta u =f(u,\partial_t u), &  x\in \mathbb{R}^n, \ t>0,\\
u(0,x)=\varepsilon u_0(x), & x\in \mathbb{R}^n, \\ \partial_t u(0,x)=\varepsilon u_1(x), & x\in \mathbb{R}^n,
\end{cases}
\end{align} has been widely investigated for several nonlinear terms $f=f(u,\partial_t u)$ over the last decades.

 For the power nonlinearity $f(u,\partial_t u)=|u|^p$, $p>1$ the small data critical exponent is the so-called Strauss exponent $p_{\Str}(n)$, that is, the positive root of the quadratic equation $(n-1)p^2-(n+1)p-2=0$. This exponent was named in  \cite{Str81} after the author's conjecture. 
 Strauss' conjecture has been proved by several authors (we address the reader to \cite[Section 20.1]{ER} 
 for  details. 
 Here, by critical exponent we mean that for $1<p\leqslant p_{\Str}(n)$ a blow-up result holds for local in time solutions under suitable sign assumptions on the data and regardless of their size, while for $p>p_{\Str}(n)$ the global in time existence of small data solutions holds in suitable function spaces. 
 
 On the other hand, for the defocusing nonlinearity $f(u,\partial_t u)=-u|u|^{p-1}$, $p>1$ the critical exponent 
  is the Sobolev exponent $p_{\Sob}(n)\doteq \frac{n+2}{n-2}$. In this case, by critical exponent we mean the existence of global in time solutions for $p<p_{\Sob}(n)$ without any assumption on the size of $\varepsilon>0$ (in the literature this kind of solutions are called large data solutions). 
  We refer to the introduction of \cite{Luc15} for further details.
 
 Finally, for the nonlinearity of derivative type $f(u,\partial_t u)=|\partial_t u|^p$ the small data critical exponent is the so-called \emph{Glassey exponent}  $$p_{\Gla}(n)\doteq \frac{n+1}{n-1}.$$ This exponent coincides with the weak solutions blow-up exponent determined in \cite{Kato80} for $f(u,\partial_t u)=|u|^p$. Coming back to $f(u,\partial_t u)=|\partial_t u|^p$, we refer to \cite{HWY12} and references therein for a summary of the known results. 
Up to our best knowledge the global existence in the supercritical case for the not radial symmetric case in high dimensions is still open. We point out that the Glassey exponent appears also in the study of models somehow related to the semilinear wave equation with nonlinearity of derivative type. For example, in \cite{LT18Glass} a blow-up result for $1<p\leqslant p_{\Gla}(n)$ has been proved for a semilinear damped wave model in the scattering case 
(see \cite{PalTak19dt} for the generalization to the case of a weakly coupled system). Moreover, in \cite{ChenPal20} the nonexistence of globally in time solutions is proved for the semilinear Moore-Gibson-Thompson equation 
 in the conservative case  with the same  assumptions on data, kind of nonlinearity and range for the exponent. 

One can try to study existence, blow up and critical exponents for weakly semilinear hyperbolic equations. A first step in this direction was made in \cite{DAn95} with $u_{tt}-a^2(t)\Delta u= -u|u|^{p-1}$, where $a(t)$ may vanish. In particular, this model includes 
 the semilinear Cauchy problem for the \emph{generalized Tricomi operator} $\mathcal{T}_\ell\doteq \partial_t^2-t^{2\ell} \Delta$, $\ell>0$, namely
\begin{align}\label{semilinear tric f}
\begin{cases} \partial_t^2 u-t^{2\ell}\Delta u =f(u,\partial_t u), &  x\in \mathbb{R}^n, \ t>0,\\
u(0,x)=\varepsilon u_0(x), & x\in \mathbb{R}^n, \\ \partial_t u(0,x)=\varepsilon u_1(x), & x\in \mathbb{R}^n.
\end{cases}
\end{align} 
Dealing with the above discussed three different kinds of nonlinearity, we expect that the corresponding critical exponents depend on $\ell$. 

For $f(u,\partial_t u)=|u|^p$ 
 the semilinear generalized Tricomi equation with power nonlinearity has been studied in several papers over the last years. Although the global existence of small data solution in the supercritical case has been proved only for $\ell=1$ and in space dimension $n=1,2$ (cf. \cite{HWY17d1,HWY18}), due to the blow-up results both in the subcritical case and in the critical case from \cite{HWY17,LinTu19}, it seems reasonable that  the critical exponent for \eqref{semilinear tric f} 
should be given by the greatest root of the quadratic equation
\begin{align*}
((\ell+1)n-1)p^2-((\ell+1)n+1-2\ell)p-2(\ell+1)=0.
\end{align*} Note that for $\ell=0$ the previous quadratic equation provides the Strauss exponent $p_{\Str}(n)$, so that this exponent is a natural generalization of the Strauss exponent.

Moreover,
in \cite{DL13} the nonexistence of global in time solutions (under suitable sign conditions) is proved by using a \emph{test function type} approach for a smaller range for $p$, namely, $1<p\leqslant p_{\Gla}((\ell+1)n)$. 

On the other hand, many papers have been devoted to the study of semilinear generalized Tricomi equation with large data and solutions defocusing nonlinearity, that is, \eqref{semilinear tric f} for $f(u,\partial_t u)=-u|u|^{p-1}$. Also in this case, a scaling of the critical exponent appears: up to dimension $n\leqslant 4$ the global existence has been proved for $p<p_{\Sob}((\ell+1)n)$. For an overview on these results we quote the most recent work  \cite{Luc19} and references therein. 


The purpose of this paper is to investigate \eqref{semilinear tric f} with nonlinearity of derivative type $f(u,\partial_t u)=|\partial_t u|^p$, which seems not present in the literature currently. According to previous results, we would expect the critical exponent to be $p_{\mathrm{Gla}}\big((\ell+1)n\big)$.  In the present paper, we prove the blow-up in finite time of a local in time solution to \eqref{semilinear CP derivative type} under suitable sign assumptions for the Cauchy data when the exponent of the nonlinearity $|\partial_t u|^p$ satisfies 
\begin{align}\label{blow-up range p}
1<p\leqslant p_{\mathrm{Gla}}\big((\ell+1)n\big)= \frac{(\ell+1)n+1}{(\ell+1)n-1}.
\end{align}
 As byproduct of the comparison argument that will be employed to prove the blow-up result, we find an upper bound estimate for the lifespan in terms of $\varepsilon$.

Let us provide an explanation on the consistency and on the reasonableness of $p_{\mathrm{Gla}}\big((\ell+1)n\big)$ as critical exponent for the semilinear Cauchy problem \eqref{semilinear CP derivative type}. 
 Although 
  a global existence result for small data solutions in the supercritical case $p>p_{\mathrm{Gla}}\big((\ell+1)n\big)$, should be considered in order to prove that this exponent is actually sharp, there are some hints that would suggest the likelihood of our conjecture. This exponent is consistent with the result for the case of the wave equation ($\ell=0$). 
  Moreover, we might interpret the parameter for the Glassey type exponent in a significant way: if we denote the quasi-homogeneous dimension of the generalized Tricomi operator $\partial_t^2-t^{2\ell}\Delta$ by $\mathcal{Q}=\mathcal{Q}(\ell)\doteq (\ell+1)n+1$ (cf. \cite{DL03,DL13,Luc18}), then, the previous exponent may be rewritten as $p_{\mathrm{Gla}}\big(\mathcal{Q}(\ell)-1\big)$. Note that, even though we are working with a nonlinearity of derivative type rather than with a power nonlinearity, this kind of exponent can be included in the class of Fujita-type critical exponents (cf. \cite[Section 2]{Luc18}).

Let us state now the main result.

\begin{theorem} \label{Thm blow-up derivative type nonlinearity} Let $n\geqslant 1$ and $\ell>0$. We assume that 
$(u_0,u_1)\in \mathcal{C}^2_0(\mathbb{R}^n)\times  \mathcal{C}^1_0(\mathbb{R}^n)$ are
 nonnegative and compactly supported in $B_R\doteq \{x\in\mathbb{R}^n: |x|<R\}$ functions.
 Let us assume that the exponent of the nonlinearity of derive type $p$ satisfies \eqref{blow-up range p}. 
 Then, there exists $\varepsilon_0=\varepsilon_0(n,p,\ell,u_0,u_1,R)>0$ such that for any $\varepsilon\in (0,\varepsilon_0]$ if $u\in \mathcal{C}^2([0,T)\times \mathbb{R}^n)$ is a local in time solution to \eqref{semilinear CP derivative type} 
  and $T=T(\varepsilon)$ is the lifespan of $u$, then, 
$u$ blows up in finite time.

 Furthermore, the following upper bound estimate for the lifespan holds
\begin{align} \label{upper bound lifespan semilinear eq}
T(\varepsilon) \leqslant \begin{cases} C \varepsilon^{-\left(\frac{1}{p-1}-\frac{(\ell+1)n-1}{2}\right)^{-1}} & \mbox{if} \ \ 1<p<p_{\Gla}\big((\ell+1)n\big), \\ \exp\big(C \varepsilon^{-(p-1)}\big) & \mbox{if} \ \ p=p_{\Gla}\big((\ell+1)n\big), \end{cases}
\end{align} where the positive constant $C$ is independent of $\varepsilon$.
\end{theorem}

We will apply a generalization of Zhou's method (see \cite{Zhou01}) for the proof of the analogous result for the semilinear wave equation with nonlinearity of derivative type. 

 In particular, instead of the classical d'Alembert's formula we shall employ  from the series of papers 
\cite{Yag04,Yag07n3} Yagdjian's integral representation formulas (obtained via a ``two-step Duhamel's principle'') for solutions to the Cauchy problem for Tricomi type equations. 
  We end up with a nonlinear ordinary integral inequality for a suitable functional related to a local solution of \eqref{semilinear CP derivative type}. Then, a comparison argument suffices to prove Theorem \ref{Thm blow-up derivative type nonlinearity}. A similar method has been applied very recently in \cite{PT19} to study the blow-up dynamic for the semilinear wave equation with time-dependent scale-invariant coefficients for the damping and mass and with nonlinearity of derivative type.

\begin{remark} The result obtained in this work can be improved in two directions. On the one hand, the validity of a global existence result for small data solutions in the supercritical case $p>p_{\mathrm{Gla}}\big(\mathcal{Q}(\ell)-1\big)$ should be proved. 

Secondly, the approach that we employed for the proof of Theorem \ref{Thm blow-up derivative type nonlinearity} strongly relies on the assumption that the Cauchy data are compactly supported. Thus, the subcritical case with not compactly supported data is open as well.
\end{remark}

\section{Fundamental tools}




\subsection{Integral representation formula}

In this section, we firstly recall an integral representation formula for the solution of the linear Cauchy problem for the generalized Tricomi equation in the one-dimensional case, namely,
\begin{align}\label{inhomog lin CP n=1}
\begin{cases} \partial_{t}^2u - t^{2\ell} \partial_{x}^2u=g(t,x), &  x\in \mathbb{R}, \ t>0,\\
u(0,x)=u_0(x), & x\in \mathbb{R}, \\ \partial_t u(0,x)=u_1(x), & x\in \mathbb{R},
\end{cases}
\end{align} where $\ell$ is a positive constant. For the proof of the representation formula one can see \cite[Theorem 3.1]{Yag04} when the data are identically 0 and \cite{Smi78} in the sourceless case. Furthermore, by following the main steps from \cite{YagGal09,PalRei18,Pal19RF} it is possible to derive the representation formula for the homogeneous case as a consequence of Yagdjian's integral formula in the inhomogeneous case with vanishing initial data.


\begin{proposition} \label{Prop representation formula 1d case}
Let $n=1$ and $\ell>0$. Let us assume $u_0\in \mathcal{C}^2_0(\mathbb{R})$, $u_1\in \mathcal{C}^1_0(\mathbb{R})$ and $g\in \mathcal{C}([0,\infty),\mathcal{C}^1( \mathbb{R}))$. 
 Then, a representation formula for the solution of \eqref{inhomog lin CP n=1} is given by
\begin{align}
u(t,x) &= a_\ell \, \phi_\ell(t)^{1-2\gamma}\int_{x-\phi_\ell(t)}^{x+\phi_\ell(t)}  u_0(y) \left(\phi_\ell(t)^2-(y-x)^2\right)^{\gamma-1}\, \mathrm{d}y \notag \\ & \qquad + b_\ell \int_{x-\phi_\ell(t)}^{x+\phi_\ell(t)}  u_1(y) \left(\phi_\ell(t)^2-(y-x)^2\right)^{-\gamma}\, \mathrm{d}y \notag \\  & \qquad + c_\ell \int_0^t \int_{x-\phi_\ell(t)+\phi_\ell(b)}^{x+\phi_\ell(t)-\phi_\ell(b)} g(b,y) E(t,x;b,y;\ell) \, \mathrm{d}y\, \mathrm{d}b, \label{representation formula 1d case}
\end{align} where the parameter $\gamma$ and the multiplicative constants $a_\ell,b_\ell, c_\ell$ are given by
\begin{align} \label{def gamma and a,b,c ell}
\gamma\doteq \frac{\ell}{2(\ell+1)}, \quad a_\ell\doteq 2^{1-2\gamma} \frac{\Gamma(2\gamma)}{\Gamma^2(\gamma)}, \quad b_\ell\doteq 2^{2\gamma-1} (\ell+1)^{1-2\gamma} \frac{\Gamma(2-2\gamma)}{\Gamma^2(1-\gamma)} , \quad  c_\ell\doteq 2^{2\gamma-1}(\ell+1)^{-2\gamma},
\end{align}
 the distance function $\phi_\ell$ is 
\begin{align} \label{def phi ell function}
 \phi_\ell(\tau)\doteq \int_0^\tau s^\ell \mathrm{d}s= \frac{\tau^{\ell+1}}{\ell+1},
\end{align} 
 and the kernel function is defined by
\begin{align}
E(t,x;b,y;\ell) & \doteq \left((\phi_\ell(t)+\phi_\ell(b))^2-(y-x)^2\right)^{-\gamma} \mathsf{F}\left(\gamma,\gamma;1; \frac{(\phi_\ell(t)-\phi_\ell(b))^2-(y-x)^2}{(\phi_\ell(t)+\phi_\ell(b))^2-(y-x)^2} \right). \label{def E(t,x;b,y)}
\end{align} Here $\Gamma(z)$ and $\mathsf{F}(a,b;c; z)$ denote the gamma function and  Gauss hypergeometric function, respectively.
\end{proposition}

In the next section, we will need to estimate from below the kernel function $E(t,x;b,y;\ell)$.

\begin{lemma} \label{Remark lower bound Hyper Gauss funct} 
Let $a\in \mathbb{R}$ and $c>0$. Then,
\begin{align*} 
\mathsf{F}(a,a;c;z)\geqslant 1  \qquad \mbox{for any }  \  z\in [0,1).
\end{align*} 
\end{lemma}

The proof of the previous estimate is based on the series expansion for $\mathsf{F}(a,a;c;z)$ (see for example \cite[Chapter 15]{OLBC10}).
By Lemma \ref{Remark lower bound Hyper Gauss funct}, it follows the lower bound estimates
\begin{align} \label{lower bound Hyper Gauss funct}
E(t,x;b,y;\ell) \geqslant \left((\phi_\ell(t)+\phi_\ell(b))^2-(y-x)^2\right)^{-\gamma}
\end{align} for any $b\in [0,t]$ and any $ y\in [x-\phi_\ell(t)+\phi_\ell(b),x+\phi_\ell(t)-\phi_\ell(b)]$.

\subsection{The curved light cone for the Tricomi equation}

 From Proposition \ref{Prop representation formula 1d case} we see that  if $\supp u_0, \supp u_1 \subset B_R$ and $\supp g\subset \{(t,x)\in [0,\infty)\times \mathbb{R}^n: |x|\leqslant \phi_\ell(t)+R\}$, then, $\supp u(t,\cdot) \subset B_{R+\phi_\ell(t)}$ for any $t\geqslant 0$.

Hence,  the forward light cone associated to the generalized Tricomi operator in the one dimensional case is 
$\{(t,x)\in [0,\infty)\times \mathbb{R}:  |x|=\phi_\ell(t)\}$. This property holds still true in higher dimensions, see for example 
 \cite{Yag04} where the representation formulae in higher dimensions are provided. 


Following the same steps as in \cite[Section 2.1]{DAn01}, it is possible to prove a local in time existence result for \eqref{semilinear CP derivative type}, regardless the size of the Cauchy data. 

Carrying out the chance of variables $v(t,x)=u(\psi(t),x)$ the semilinear equation in \eqref{semilinear CP derivative type} can be rewritten as
\begin{align*}
v_{tt}= (\psi)^{2\ell}(\psi')^2 \Delta v+\frac{\psi''}{\psi'} v_t+(\psi')^{2-p}|v_t|^p.
\end{align*} Therefore, choosing $\psi$ so that $\psi'=\psi^{-\ell}$, that is, $\psi$ is the inverse function of $\phi_\ell$, we have that $v$ solves 
\begin{align*} 
v_{tt}-\Delta v+ \frac{\mu_\ell}{t} \,\partial_t v = c_{\ell,p} t^{\mu_\ell (p-2)} |\partial_t v|^p,
\end{align*} where $\mu_\ell\doteq \frac{\ell}{\ell+1}$ and $c_{\ell,p}\doteq (\ell+1)^{\mu_\ell(p-2)}$,  which is a semilinear Euler-Darboux-Poisson equation. 

For $t>0$ the equation is strictly hyperbolic and the classical theory applies. In particular, $v=v(t,x)$ has  finite speed of propagation property. Furthermore, the light cone for $v$ is given by $\{(t,x)\in [0,\infty)\times \mathbb{R}^n: |x|= t\}$, so, coming back to $u=u(t,x)$, if $\supp u_0, \supp u_1 \subset B_R$  then, 
\begin{align} \label{supp cond}\supp u(t,\cdot) \subset B_{R+\phi_\ell(t)} \quad \mbox{for any} \ \  t\geqslant 0.
\end{align} 


\begin{remark}
According to the previously introduced change of variables, applying the techniques from our approach, we may also obtain a blow-up result for the semilinear Euler-Darboux-Poisson equation with nonlinearity of derivative type. 
\end{remark}

\section{Proof of Theorem \ref{Thm blow-up derivative type nonlinearity}} 

Let $u$ be a local (in time) solution to the Cauchy problem \eqref{semilinear CP derivative type}. We introduce an auxiliary function which depends on the time variable and only on the first space variable, by integrating $u$ with respect to the remaining $(n-1)$ spatial variables. This means that, if we denote $x=(z,w)$ with $z\in \mathbb{R}$ and $w\in \mathbb{R}^{n-1}$, then, we deal with the function
\begin{align*}
\mathcal{U}(t,z)\doteq \int_{\mathbb{R}^{n-1}} u(t,z,w) \, \mathrm{d} w \qquad \mbox{for any} \ \  t>0, z\in\mathbb{R}.
\end{align*} 
 Hereafter, we will deal formally only with the case $n\geqslant 2$ for the sake of brevity, although one can proceed exactly in the same way for $n=1$ by working with $u$ in place of $\mathcal{U}$. In order to describe the initial values of $\mathcal{U}$, we introduce
\begin{align*}
\mathcal{U}_0(z)\doteq \int_{\mathbb{R}^{n-1}} u_0(z,w) \, \mathrm{d} w , \quad \mathcal{U}_1(z)\doteq \int_{\mathbb{R}^{n-1}} u_1(z,w) \, \mathrm{d} w\qquad \mbox{for any} \ \  z\in\mathbb{R}.
\end{align*} Since we assume that $u_0,u_1$ are compactly supported with support contained in $B_R$, it follows that $\mathcal{U}_0,\mathcal{U}_1$ are compactly supported in $(-R,R)$. Similarly, 
due to the property of finite speed of propagation of perturbations for $u$, from \eqref{supp cond} we have
\begin{align} \label{supp mathcal U}
\supp \mathcal{U}(t,\cdot)\subset (-(R+\phi_\ell(t)),R+\phi_\ell(t)) \qquad \mbox{for any} \ \  t>0.
\end{align}  
Consequently, $\mathcal{U}$ solves the following Cauchy problem
\begin{align*} 
\begin{cases}
\partial_t^2\mathcal{U}- t^{2\ell}\partial_z^2\mathcal{U}=\int_{\mathbb{R}^{n-1}} |\partial_t u(t,z,w)|^p \, \mathrm{d} w , &  z\in \mathbb{R} , \ t>0,\\
\mathcal{U}(0,z)= \varepsilon\, \mathcal{U}_0(z), & z\in \mathbb{R} , \\ \partial_t\mathcal{U}(0,z)= \varepsilon\, \mathcal{U}_1(z) , & z\in \mathbb{R}.
\end{cases}
\end{align*} By Proposition \ref{Prop representation formula 1d case} it follows
\begin{align*}
\mathcal{U}(t,z) &= a_\ell \, \varepsilon \, \phi_\ell(t)^{1-2\gamma}\int_{z-\phi_\ell(t)}^{z+\phi_\ell(t)}  \mathcal{U}_0(y) \left(\phi_\ell(t)^2-(y-z)^2\right)^{\gamma-1} \mathrm{d}y \notag \\ & \qquad + b_\ell \, \varepsilon \int_{z-\phi_\ell(t)}^{z+\phi_\ell(t)}  \mathcal{U}_1(y) \left(\phi_\ell(t)^2-(y-z)^2\right)^{-\gamma} \mathrm{d}y \notag \\  & \qquad + c_\ell \int_0^t \int_{z-\phi_\ell(t)+\phi_\ell(b)}^{z+\phi_\ell(t)-\phi_\ell(b)} \int_{\mathbb{R}^{n-1}} |\partial_t u(b,y,w)|^p \, \mathrm{d} w \, E(t,z;b,y;\ell) \, \mathrm{d}y\, \mathrm{d}b,
\end{align*} where the kernel function $E$ is defined by \eqref{def E(t,x;b,y)}.

Due to the sign assumption for $u_0,u_1$ it follows that $\mathcal{U}_0,\mathcal{U}_1$ are nonnegative functions. Therefore, using the fact that $\gamma\in (0,1)$ and estimating the kernel functions in the integrals containing $\mathcal{U}_0$ and $\mathcal{U}_1$ from below, we obtain
\begin{align}
\mathcal{U}(t,z) &\geqslant 2^{\gamma-1} a_\ell \, \varepsilon \, \phi_\ell(t)^{-\gamma}\int_{z-\phi_\ell(t)}^{z+\phi_\ell(t)}  \mathcal{U}_0(y) \big(\phi_\ell(t)-z+y\big)^{\gamma-1} \, \mathrm{d}y  \notag \\
  & \qquad  +2^{-\gamma}  b_\ell \, \varepsilon  \, \phi_\ell(t)^{-\gamma} \int_{z-\phi_\ell(t)}^{z+\phi_\ell(t)}  \mathcal{U}_1(y) \big(\phi_\ell(t)-z+y\big)^{-\gamma} \, \mathrm{d}y \notag \\  & \qquad + c_\ell \int_0^t \int_{z-\phi_\ell(t)+\phi_\ell(b)}^{z+\phi_\ell(t)-\phi_\ell(b)} \int_{\mathbb{R}^{n-1}} |\partial_t u(b,y,w)|^p \, \mathrm{d} w \, E(t,z;b,y;\ell) \, \mathrm{d}y\, \mathrm{d}b. \label{1st lb mathcal U}
\end{align} 
Let us investigate the behavior of the terms
\begin{align*}
 J(t,z) & \doteq   \phi_\ell(t)^{-\gamma}\int_{z-\phi_\ell(t)}^{z+\phi_\ell(t)}  \left( 2^{\gamma-1} a_\ell \mathcal{U}_0(y) \big(\phi_\ell(t)-z+y\big)^{\gamma-1}  +2^{-\gamma} b_\ell  \mathcal{U}_1(y) \big(\phi_\ell(t)-z+y\big)^{-\gamma}\right)  \mathrm{d}y , \\
 I(t,z) & \doteq  c_\ell \int_0^t \int_{z-\phi_\ell(t)+\phi_\ell(b)}^{z+\phi_\ell(t)-\phi_\ell(b)} \int_{\mathbb{R}^{n-1}} |\partial_t u(b,y,w)|^p \, \mathrm{d} w \, E(t,z;b,y;\ell) \, \mathrm{d}y\, \mathrm{d}b.
\end{align*}
On the characteristic line $\phi_\ell(t)-z=R$ and for $z\geqslant R$, it holds $[-R,R]\subset [z-\phi_\ell(t),z+\phi_\ell(t)]$. Since $\supp \mathcal{U}_0,\supp \mathcal{U}_1\subset (-R,R)$, we may estimate
\begin{align*}
 J(t,z) &  \geqslant \phi_\ell(t)^{-\gamma}\int_{-R}^{R}  \left( 2^{2(\gamma-1)} a_\ell R^{\gamma-1} \mathcal{U}_0(y)  + 2^{-2\gamma} b_\ell R^{-\gamma}  \mathcal{U}_1(y) \right)  \mathrm{d}y \\ 
  & =  (z+R)^{-\gamma}\int_{\mathbb{R}}  \big( 2^{2(\gamma-1)} a_\ell R^{\gamma-1} \mathcal{U}_0(y)  +  2^{-2\gamma} b_\ell R^{-\gamma}  \mathcal{U}_1(y)\big) \, \mathrm{d}y \\ 
  &\geqslant K  (z+R)^{-\gamma}\int_{\mathbb{R}^n}  \big(  u_0(y) + u_1(y)\big) \, \mathrm{d}y,
\end{align*}  where $K=K(R,\ell)=\min\{2^{2(\gamma-1)} a_\ell R^{\gamma-1},2^{-2\gamma} b_\ell R^{-\gamma}\}$.

Next we estimate the term $I(t,z)$. Due to the support condition  $\supp \partial_t u(t,\cdot) \subset B_{R+\phi_\ell(t)}$ it results $$\supp \partial_t u(t,z,\cdot)\subset \big\{w\in \mathbb{R}^{n-1}: |w|\leqslant \left((R+\phi_\ell(t))^2-z^2\right)^{1/2} \big\} \quad \mbox{for any} \ \ t>0, z\in \mathbb{R}.$$ Then, by H\"older's inequality, we have
\begin{align*}
  |\partial_t \mathcal{U} (b,y)| & =\Big| \int_{\mathbb{R}^{n-1}} \partial_t u(b, y,w) \, \mathrm{d}w\Big| \\ &  \leqslant \bigg(\int_{\mathbb{R}^{n-1}} |\partial_t u(b, y,w)|^p \, \mathrm{d}w\bigg)^{\frac{1}{p}} \big(\meas_{n-1}\big(\supp \partial_t u(b,y,\cdot)\big)\big)^{1-\frac{1}{p}} \\ & \lesssim \left((R+\phi_\ell(b))^2-y^2\right)^{\frac{n-1}{2}\left(1-\frac{1}{p}\right)}\bigg(\int_{\mathbb{R}^{n-1}} |\partial_t u(b, y,w)|^p \, \mathrm{d}w\bigg)^{\frac{1}{p}} .
\end{align*}
Henceforth, the unexpressed multiplicative constants will depend on $n,\ell,R,p$. 
 Hence,
\begin{align*}
\int_{\mathbb{R}^{n-1}} |\partial_t u(b, y,w)|^p \, \mathrm{d}w  & \gtrsim  \left((R+\phi_\ell(b))^2-y^2\right)^{-\frac{n-1}{2}(p-1)} |\partial_t \mathcal{U} (b,y)|^p,
\end{align*}  which implies in turn
\begin{align*}
I(t,z) &  \gtrsim  \int_0^t \int_{z-\phi_\ell(t)+\phi_\ell(b)}^{z+\phi_\ell(t)-\phi_\ell(b)} \left((R+\phi_\ell(b))^2-y^2\right)^{-\frac{n-1}{2}(p-1)} |\partial_t \mathcal{U} (b,y)|^p  E(t,z;b,y;\ell) \, \mathrm{d}y\, \mathrm{d}b \\ 
& =  \int_{z-\phi_\ell(t)}^{z+\phi_\ell(t)} \int_0^{\phi^{-1}_\ell(\phi_\ell(t)-|z-y|)} \left((R+\phi_\ell(b))^2-y^2\right)^{-\frac{n-1}{2}(p-1)} |\partial_t \mathcal{U}(b,y)|^p  E(t,z;b,y;\ell)  \, \mathrm{d}b \, \mathrm{d}y ,
\end{align*} where we used Fubini's theorem in the last equality. 

Hereafter, we work on the characteristic line $\phi_\ell(t)-z=R$ and for $z\geqslant R$. Thus, shrinking the domain of integration, we get
\begin{align*}
I(\phi_\ell^{-1}(z+R),z) &  \gtrsim \int_{R}^{z} \int_{\phi_\ell^{-1}(y-R)}^{\phi_\ell^{-1}(y+R)} \left((R+\phi_\ell(b))^2-y^2\right)^{-\frac{n-1}{2}(p-1)} |\partial_t \mathcal{U}(b,y)|^p  E(t,z;b,y;\ell)  \, \mathrm{d}b \, \mathrm{d}y \\
& \gtrsim \int_{R}^{z} \left(R+y\right)^{-\frac{n-1}{2}(p-1)} \int_{\phi_\ell^{-1}(y-R)}^{\phi_\ell^{-1}(y+R)} |\partial_t \mathcal{U}(b,y)|^p  E(t,z;b,y;\ell)  \, \mathrm{d}b \, \mathrm{d}y.
\end{align*}  

Let us fix $(t,z;y,b)$ such that $\phi_\ell(t)-z=R$ and $z\geqslant R$, $y\in [R,z]$, $ b\in [\phi_\ell^{-1}(y-R),\phi_\ell^{-1}(y+R)]$, then, from \eqref{lower bound Hyper Gauss funct} we get
\begin{align*}
E(t,z;b,y;\ell) & \geqslant  \left((\phi_\ell(t)+\phi_\ell(b))^2-(z-y)^2\right)^{-\gamma}  = \left((z+R+\phi_\ell(b))^2-(z-y)^2\right)^{-\gamma}  \\
 & \geqslant  \left((z+y+2R)^2-(z-y)^2\right)^{-\gamma} = 2^{-2\gamma} (y+R)^{-\gamma} (z+R)^{-\gamma}.
\end{align*}
 Consequently, 
 we obtain
\begin{align} \label{lower bound I}
I(\phi_\ell^{-1}(z+R),z)   
& \gtrsim (z+R)^{-\gamma} \int_{R}^{z} \left(R+y\right)^{-\frac{n-1}{2}(p-1)-\gamma} \int_{\phi_\ell^{-1}(y-R)}^{\phi_\ell^{-1}(y+R)} |\mathcal{U}_t(b,y)|^p   \, \mathrm{d}b \, \mathrm{d}y.
\end{align} 
 Applying Jensen's inequality and $\phi_\ell^{-1}(\tau)=\big((\ell+1)\tau\big)^{\frac{1}{\ell+1}}$, we find
 \begin{align}
\left| \int_{\phi_\ell^{-1}(y-R)}^{\phi_\ell^{-1}(y+R)} \partial_t \mathcal{U}(b,y)   \, \mathrm{d}b \,\right|^p & \leqslant \big(\phi_\ell^{-1}(y+R)-\phi_\ell^{-1}(y-R)\big)^{p-1} \int_{\phi_\ell^{-1}(y-R)}^{\phi_\ell^{-1}(y+R)} |\mathcal{U}_t(b,y)|^p   \, \mathrm{d}b \notag \\
& \leqslant (2R(\ell+1))^{\frac{p-1}{\ell+1}} \int_{\phi_\ell^{-1}(y-R)}^{\phi_\ell^{-1}(y+R)} |\partial_t \mathcal{U}(b,y)|^p   \, \mathrm{d}b. \label{Jensen mathcal U_t}
 \end{align}
 Combining \eqref{lower bound I}, \eqref{Jensen mathcal U_t} and the fundamental theorem of calculus, we arrive at
\begin{align*}
I(\phi_\ell^{-1}(z+R),z)   
& \gtrsim (z+R)^{-\gamma} \int_{R}^{z} \left(R+y\right)^{-\frac{n-1}{2}(p-1)-\gamma} \left| \int_{\phi_\ell^{-1}(y-R)}^{\phi_\ell^{-1}(y+R)} \partial_t \mathcal{U}(b,y)   \, \mathrm{d}b \,\right|^p \mathrm{d}y \\
& = (z+R)^{-\gamma} \int_{R}^{z} \left(R+y\right)^{-\frac{n-1}{2}(p-1)-\gamma} |\mathcal{U}(\phi_\ell^{-1}(y+R),y) |^p \, \mathrm{d}y, 
\end{align*} where in the second step we used $\mathcal{U}(\phi_\ell^{-1}(y-R),y)=0$ due to \eqref{supp mathcal U}.
  Combining the lower bound estimates for $J$ and $I$, on the characteristic $\phi_\ell(t)-z=R$ and for $z\geqslant R$, 
  it results
\begin{align} \label{fundamental inequality for mathcal U}
(R+z)^{\gamma}\mathcal{U}(\phi_\ell^{-1}(z+R),z) & \geqslant K \varepsilon \, \|  u_0 +   u_1 \|_{L^1(\mathbb{R}^n)}  + C \int_{R}^{z} \left(R+y\right)^{-\frac{n-1}{2}(p-1)-\gamma} |\mathcal{U}(\phi_\ell^{-1}(y+R),y) |^p \, \mathrm{d}y,
\end{align} where $C=C(n,\ell,R,p)>0$ is the unexpressed multiplicative constant appearing in the estimate from below of $I(\phi_\ell^{-1}(z+R),z)$. 
We define $$U(z)\doteq (R+z)^{\gamma} \mathcal{U}(\phi_\ell^{-1}(z+R),z).$$ We shall use the dynamic of this function to prove the blow-up result. We may rewrite  \eqref{fundamental inequality for mathcal U} as
\begin{align}\label{fundamental inequality for U}
U(z)\geqslant K \varepsilon \, \|  u_0 +  u_1 \|_{L^1(\mathbb{R}^n)}  + C \int_{R}^{z} \left(R+y\right)^{-\frac{n-1}{2}(p-1)-\gamma(p+1)} |U(y) |^p \, \mathrm{d}y \qquad \mbox{for}  \ z\geqslant R.
\end{align} 
Now we apply a comparison argument to $U$. We introduce the auxiliary function $G$ as follows:
\begin{align*}
G(z)\doteq   M \varepsilon + C \int_{R}^{z} \left(R+y\right)^{-\frac{n-1}{2}(p-1)-\gamma(p+1)} |U(y) |^p \, \mathrm{d}y \qquad \mbox{for}  \ z\geqslant R,
\end{align*} where $M\doteq K \|  u_0 +  u_1 \|_{L^1(\mathbb{R}^n)}$. Then, by \eqref{fundamental inequality for U} we get $U\geqslant G$. Moreover, $G$ satisfies the differential inequality
\begin{align*}
G'(z) & = C \left(R+z\right)^{-\frac{n-1}{2}(p-1)-\gamma(p+1)} |U(z) |^p \\
& \geqslant C \left(R+z\right)^{-\frac{n-1}{2}(p-1)-\gamma(p+1)} (G(z) )^p . 
\end{align*} 
Let us consider the initial value $G(R)=M\varepsilon$.

If 
\begin{align} \label{critical expression}
-\tfrac{n-1}{2}(p-1)-\gamma(p+1)=-1,
\end{align} then, 
\begin{align*}
(M\varepsilon)^{1-p}-G(z)^{1-p}\geqslant C (p-1)\log \left(\frac{R+z}{2R}\right) .
\end{align*} Otherwise, being $G$ a positive function, we have 
\begin{align} \label{subcritical expression}
(M\varepsilon)^{1-p}-G(z)^{1-p} \geqslant  \frac{C(\ell+1)}{\frac{1}{p-1}-\frac{(\ell+1)n-1}{2}} \Big((R+z)^{1-2\gamma-\frac{n+2\gamma-1}{2}(p-1)}-(2R)^{1-2\gamma-\frac{
n+2\gamma-1}{2}(p-1)}\Big).
\end{align}

We point out that \eqref{critical expression} is equivalent to require $p=p_{\Gla}\big((\ell+1)n\big) $. Hence, if $p\in \big(1,p_{\Gla}\big((\ell+1)n\big)\big)$ the multiplicative factor on the right-hand side of \eqref{subcritical expression} is positive and we can choose $\varepsilon_0=\varepsilon_0(n,p,\ell,u_0,u_1,R)$ sufficiently small such that for $\varepsilon\in (0,\varepsilon_0]$ we obtain
\begin{align}
G(z) & \geqslant \left[(M\varepsilon)^{1-p} +\frac{C(\ell+1)}{\frac{1}{p-1}-\frac{(\ell+1)n-1}{2}} \Big((2R)^{1-2\gamma-\frac{n+2\gamma-1}{2}(p-1)}-(R+z)^{1-2\gamma-\frac{n+2\gamma-1}{2}(p-1)}\Big) \right]^{-\frac{1}{p-1}}  \notag \\
&\geqslant \left[2(M\varepsilon)^{1-p}- \frac{C(\ell+1)}{\frac{1}{p-1}-\frac{(\ell+1)n-1}{2}}(R+z)^{1-2\gamma-\frac{n+2\gamma-1}{2}(p-1)} \right]^{-\frac{1}{p-1}}. \label{estimate G below subcrit}
\end{align} 

In the critical case $p=p_{\Gla}\big((\ell+1)n\big)$, we have that
 \begin{align*}
U(z)\geqslant G(z) & \geqslant  \left[(M\varepsilon)^{1-p}-C(p-1)\log \left(\tfrac{R+z}{2R}\right)\right]^{-\frac{1}{p-1}}
 \end{align*} implies the blow-up in finite time of $U(z)$  and the lifespan estimate $$T(\varepsilon)\lesssim \exp\left(\widetilde{C}\varepsilon^{-(p-1)}\right).$$
On the other hand, in the subcritical case $p\in \big(1,p_{\Gla}\big((\ell+1)n\big)\big)$, in \eqref{estimate G below subcrit} the right-hand side blows up for
$$\phi_\ell(t)=R+z\simeq \varepsilon^{-(\ell+1)\left(\frac{1}{p-1}-\frac{(\ell+1)n-1}{2}\right)^{-1}},$$ where we used the relation
\begin{align*}
1-2\gamma-\tfrac{n+2\gamma-1}{2}(p-1) =\tfrac{1}{\ell+1}\left(1-\tfrac{(\ell+1)n-1}{2}(p-1)\right).
\end{align*}
 Being $U \geqslant G$, then, $U$ blows up in finite time and the upper bound for the lifespan is given by $$T(\varepsilon)\lesssim  \varepsilon^{-\left(\frac{1}{p-1}-\frac{(\ell+1)n-1}{2}\right)^{-1}}.$$ 
  Therefore, the proof of Theorem \ref{Thm blow-up derivative type nonlinearity} is completed.

\section*{Acknowledgments}


A. Palmieri is supported by the GNAMPA project `Problemi stazionari e di evoluzione nelle equazioni di campo nonlineari dispersive'.
S. Lucente is supported by the PRIN 2017 project `Qualitative and quantitative aspects of nonlinear PDEs' and by the GNAMPA project `Equazioni di tipo dispersivo: teoria e metodi'.


\addcontentsline{toc}{chapter}{Bibliography}


\begin{thebibliography}{00}



\bibitem{ChenPal20} Chen, W., Palmieri, A.:
\newblock{A blow – up result for the semilinear Moore-Gibson-Thompson equation with nonlinearity of derivative type in the conservative case.}
\newblock{\em Evol. Equ. Control Theory.} (2020) to appear 


\bibitem{DL13} D’Abbicco, M., Lucente, S.:
\newblock{A modified test function method for damped wave equations.} 
\newblock{\em Adv. Nonlinear Stud.} {\bf 13}(4), 867--892  (2013)




\bibitem{DL03} D'Ambrosio, L., Lucente, S.:
\newblock{Nonlinear Liouville theorems for Grushin and Tricomi operators.}
\newblock{\em J. Differential Equations} {\bf 193},  511--541 (2003)

\bibitem{DAn95} D'Ancona, P.:
\newblock{A note on a theorem of J\"orgens.}
\newblock{\em Math. Z.} {\bf 218}, 239-252  (1995)

\bibitem{DAn01} D'Ancona, P., Di Giuseppe, A.:
\newblock{Global Existence with Large Data for a Nonlinear Weakly Hyperbolic Equation.}
\newblock{\em Math. Nachr.} {\bf 231}, 5-23  (2001)




\bibitem{ER} Ebert, M.R., Reissig, M.:
\newblock{\em Methods for partial differential equations. Qualitative properties of solutions, phase space analysis, semilinear models.}
Birkha\"user/Springer, Cham, 2018










\bibitem{HWY17} He, D., Witt, I., Yin, H.: .  56, 21 (2017). 
\newblock{On the global solution problem for semilinear generalized Tricomi equations, I.}
\newblock{\em Calc. Var.} {\bf 56}, 21 (2017)

\bibitem{HWY17d1} He, D., Witt, I., Yin, H.:
\newblock{On semilinear Tricomi equations with critical exponents or in two space dimensions.}
\newblock{\em J. Differential Equations} {\bf 263}(12), 8102--8137 (2017)

\bibitem{HWY18} He, D., Witt, I., Yin, H.:
\newblock{On semilinear Tricomi equations in one space dimension.}
 \newblock{Preprint, arXiv:1810.12748, 2018}



\bibitem{HWY12}  Hidano, K., Wang,  C.,  Yokoyama, K.: 
\newblock{ The Glassey conjecture with radially symmetric data.} \newblock{\em J. Math. Pures Appl.}  {\bf 98}(9), 518--541 (2012)









\bibitem{Kato80}  Kato, T.:
\newblock{ Blow-up of solutions of some nonlinear hyperbolic equations.} \newblock{\em Comm. Pure Appl. Math.} {\bf 33}(4), 501--505 (1980)









\bibitem{LT18Glass}  Lai, N.A., Takamura, H.:
\newblock{Nonexistence of global solutions of nonlinear wave equations with weak time-dependent damping related to Glassey's conjecture.} \newblock{\em Differential and Integral Equations} {\bf 32}(1-2), 37--48 (2019) 




\bibitem{LinTu19} Lin, J., Tu, Z.:
\newblock{Lifespan of semilinear generalized Tricomi equation with Strauss type exponent.}
\newblock{Preprint, arXiv:1903.11351v2 , 2019}  





\bibitem{Luc15} Lucente, S.:
\newblock{Large data solutions for critical semilinear weakly hyperbolic equations.} In \newblock{\em Proceeding of the Conference Complex Analysis \& Dinamical System VI,  Contemporary Mathematics}  {\bf 653} 
(2015), 251--276


\bibitem{Luc18} Lucente, S.:
\newblock{Critical exponents and where to find them.}
\newblock{\em Bruno Pini Mathematical Analysis Seminar}, [S.l.],  102-114 (2018) 



\bibitem{Luc19} Lucente, S.: 
\newblock{4D Semilinear Weakly Hyperbolic Wave Equations.} In: D'Abbicco M., Ebert M., Georgiev V., Ozawa T. (eds) \newblock{\em New Tools for Nonlinear PDEs and Application.} Trends in Mathematics. Birkh\"auser, Cham (2019)






\bibitem{OLBC10} Olver, F.W.J., Lozier, D.W., Boisvert, R.F., Clark, C.W:
\newblock{\em NIST Handbook of Mathematical Functions.}
\newblock{Cambridge University Press} (2010), p. 966


	



\bibitem{Pal19RF} Palmieri, A.:
\newblock{An integral representation formula for the solutions of a wave equation with time-dependent damping and mass in the scale-invariant case.}  \newblock{Preprint, arXiv:1905.02408, 2019}


\bibitem{PalRei18} Palmieri, A.,  Reissig, M.:
\newblock{A competition between Fujita and Strauss type exponents for blow-up of semi-linear wave equations with scale-invariant damping and mass.}
\newblock{\em J. Differential Equations} {\bf 266 }, 1176--1220 (2019)



\bibitem{PalTak19dt}  Palmieri, A., Takamura, H.:
	\newblock{Nonexistence of global solutions for a weakly coupled system of semilinear damped wave equations of derivative type in the scattering case.} 
	\newblock{\em Mediterr. J. Math} {\bf 17} 13 (2020) 
	



\bibitem{PT19}  Palmieri, A.,  Tu, Z.:
\newblock{A blow-up result for a semilinear wave equation with scale-invariant damping and mass and nonlinearity of derivative type.} \newblock{Preprint, arXiv:1905.11025v2, 2019}  










\bibitem{Smi78}  Smirnov, M.M.:
\newblock{ Equations of Mixed Type.}
\newblock{ Translations of Mathematical Monographs, vol. 51, American Mathematical Society}, Providence, RI, 1978


\bibitem{Str81} Strauss, W.A.:
\newblock{Nonlinear scattering theory at low energy.}
\newblock{\em J. Funct. Anal.} {\bf 41}(1), 110--133  (1981)
	






	









\bibitem{Yag04}  Yagdjian, K.:
\newblock{ A note on the fundamental solution for the Tricomi-type equation in the hyperbolic domain.}
\newblock{\em  J. Differential Equations} {\bf 206}, 227--252 (2004)
%




\bibitem{Yag07n3} Yagdjian, K.:
\newblock{The self-similar solutions of the one-dimensional semilinear Tricomi-type equations.}
\newblock{\em J. Differential Equations} {\bf 236}, 82--115 (2007)

%
%
%
%
%
%
%
%
\bibitem{YagGal09} Yagdjian, K., Galstian, A.:
\newblock{ Fundamental Solutions for the Klein-Gordon Equation in de Sitter Spacetime,}
\newblock{\em  Comm. Math. Phys.} {\bf 285} , 293--344 (2009)




%

\bibitem{Zhou01} Zhou, Y.: \newblock{ Blow up of solutions to the Cauchy problem for nonlinear wave equations.} \newblock{\em Chin. Ann. Math. Ser. B} {\bf 22}(3), 275--280  (2001)
	

	



\end{thebibliography}
\end{document}